\documentclass{article}
\usepackage{amssymb}
\usepackage{amsmath}
\usepackage{longtable}
\usepackage{booktabs}
\usepackage[ruled]{algorithm2e}
\usepackage{multirow,bigstrut,bigdelim}

\newtheorem{Thm}{Theorem}[section]
\newtheorem{Cor}[Thm]{Corollary}
\newtheorem{Def}[Thm]{Definition}

\newtheorem{La}[Thm]{Lemma}
\newtheorem{Rem}[Thm]{Remark}
\newtheorem{hypothesis}[Thm]{Hypothesis}

\newenvironment{Prf}{\noindent\textbf{Proof.}}{\hfill $\Box$ \medskip}

\newcommand{\F}{\mathbb{F}}
\newcommand{\Z}{\mathbb{Z}}

\newcommand{\calC}{\mathcal{C}}
\newcommand{\GL}{\textrm{GL}}

\newcommand{\GU}{\textrm{GU}}

\newcommand{\SO}{\textrm{SO}}

\newcommand{\Sp}{\textrm{Sp}}

\numberwithin{equation}{section}

\begin{document}

\title{Abundant $p$-singular elements in \\ finite classical groups}
\author{Alice C. Niemeyer, Tomasz Popiel \& Cheryl E. Praeger\footnote{The University of Western Australia, 35 Stirling Highway, Crawley, WA 6009, Australia (all authors);  
King Abdulaziz University, Jeddah, Saudi Arabia (C.~E.~Praeger); \newline 
e-mail: alice.niemeyer@uwa.edu.au, tomasz.popiel@uwa.edu.au, cheryl.praeger@uwa.edu.au.}}
\maketitle

\begin{abstract}
\noindent 
In 1995, Isaacs, Kantor and Spaltenstein proved that for a finite simple classical group $G$ defined over a field with $q$ elements, and for a prime divisor $p$ of $|G|$ distinct from the characteristic, the proportion of $p$-singular elements in $G$ (elements with order divisible by $p$) is at least a constant multiple of $(1-1/p)/e$, where $e$ is the order of $q$ modulo $p$.  
Motivated by algorithmic applications, we define a subfamily of $p$-singular elements, called $p$-abundant elements, which leave invariant certain ``large'' subspaces of the natural $G$-module.
We find explicit upper and lower bounds for the proportion of $p$-abundant elements in $G$, and prove that it approaches a (positive) limiting value as the dimension of $G$ tends to infinity. 
It turns out that the limiting proportion of $p$-abundant elements is at least a constant multiple of the Isaacs--Kantor--Spaltenstein lower bound for the proportion of {\em all} $p$-singular elements. \\


\end{abstract}

\section{Introduction} \label{sec:intro}

Given a  prime $p$ dividing  the order of  a finite group $G$, what proportion of elements 
in $G$ are $p$-singular? That is, what
proportion of elements in $G$ have order divisible by $p$?  Isaacs et~al.~\cite{IKS95} 
considered this problem for permutation groups of degree $n$ and
proved that  the proportion is at  least $1/n$. At the  heart of their
proof is consideration of the case  where $G$ is a finite simple group
of Lie type, and more particularly a finite simple $d$-dimensional classical
group. In  this case  they
obtained for the proportion of $p$-singular elements in $G$ 
a  lower  bound of  the  form  $(1-1/p)c/d$  for some  constant  c,
independently of the type  of $G$ and of the size $q$  of the field over
which $G$ is defined.  A closer inspection of their proof reveals that
their lower bound can in fact be written as a constant times $(1-1/p)/e$, where 
$e$ is the order of $q$ modulo $p$ (the least positive  integer for which $p$ divides $q^e-1$);  
we thank Klaus  Lux for pointing this out to us.

Estimates for proportions of $p$-singular elements are important in 
complexity analyses of numerous algorithms in computational group theory.  In particular, 
the first and third authors' algorithm~\cite{NPClassRec} 
for recognising finite classical groups in their natural representations 
relies on finding by repeated independent random selection from $G$ 
elements with orders divisible by primes $p$ for which $e$ is greater than 
half the dimension $d$,  
and  there exists an efficient practical algorithm
for testing  whether  an   element  has  this  property~\cite{NPAlg}.   However,  these particular
$p$-singular elements  are relatively  scarce: they arise with frequency proportional 
to $1/e = O(1/d)$, whereas the work of Isaacs et~al.~\cite{IKS95} suggests 
that, in general, $p$-singular elements are more frequent when $e$ is smaller.
Moreover,  the restriction $e>d/2$  only  allows us  to  identify
elements  with orders  divisible by  {\em certain} primes, namely those primes $p$ for which the 
order of $q$ modulo $p$ is greater than $d/2$.  

These shortcomings motivated us to
seek, for {\em all} values of $e$, a class of $p$-singular elements
that  arise with  frequency proportional to $1/e$  and can  be efficiently recognised
algorithmically, with the hope being that such elements might lead to improved
recognition algorithms for finite classical groups. 
Experimental evidence gathered by the first author suggested that  a 
particular type of $p$-singular element,  which  we  term
\emph{$p$-abundant} (see Definition~\ref{def:pab-galois}), 
arises  approximately with frequency proportional to  $1/e$.  The
theoretical analysis  presented in this  paper proves that this  is indeed
the   case,   for  the   groups   in  Table~\ref{tab:G}.    Moreover, the 
$p$-abundant elements  are indeed readily identifiable  computationally; for
example, from their characteristic  polynomials.   Algorithms  for  this 
are presented in  a companion paper \cite{alg}.

Our estimates  for the proportion  of $p$-abundant elements  are very
precise:  we determine in  Theorem~\ref{Thm:intro} both  upper  and lower
bounds  and,  in  particular,  the  exact asymptotic  value  for  this
proportion.   Precision  of  this  kind   seems  to  be  rare  in  the
literature: estimates for the proportions of various kinds of elements in 
finite groups tend to focus on lower bounds, and good upper bounds are
rarely given. Our endeavour to obtain such precise bounds  drove the development
of the  methodology presented in the first and third authors' paper~\cite{NPquokka}, 
which  underlies the proofs
in this paper. This theory was in turn inspired by and developed from
the methods used by Isaacs et~al.~\cite{IKS95} and an earlier application  
by the first and third authors in collaboration with L\"ubeck~\cite{LNP}. 
We note that this method was used earlier by Lehrer~\cite{Lehrer92,Lehrer98}  
to study the  representations of  finite Lie type groups.

Our method  requires us to sum  over the lengths  of certain conjugacy
classes  of the corresponding  Weyl group,  weighted by  proportions of $p$-abundant elements in 
matching maximal tori.   Whereas  previous  applications~\cite{LNP,NPP}  
approximated  the corresponding
expressions  by replacing all  weighting factors  with a  common lower
bound,  here we have  to be much more careful with our estimates.
Our  results  highlight the  power  of  the  method of~\cite{NPquokka} to  obtain  exact
asymptotic values for element proportions, and rely on some delicate technical
lemmas. They constitute the first application of this theory achieving such precision. 

\begin{table}[!t]
\begin{center}
\begin{tabular}{lll}
\toprule
$G$ & $\delta$ & $d$ \\
\midrule
$\GL_n(q)$ & $1$ & $n$ \\ 
$\GU_n(q)$ & $2$ & $n$ \\ 
$\Sp_{2n}(q)$ & $1$ & $2n$ \\ 
$\SO_{2n+1}(q)$ & $1$ & $2n+1$ \\ 
$\SO^\pm_{2n}(q)$ & $1$ & $2n$ \\ 
\bottomrule
\end{tabular}
\end{center}
\caption{The finite classical groups considered in this paper.}\label{tab:G}
\end{table}

We now formally define $p$-abundant elements in finite classical groups and state our results concerning proportions of these elements.

\begin{Def}\label{def:pab-galois}
\textnormal{ Let  $q$ be  a prime power,  $n$ a positive  integer, and
  $G$,   $\delta$,   $d=d(n)$   as    in   one   of   the   lines   of
  Table~\ref{tab:G}.   Let   $V=V(d,q^\delta)$   denote  the   natural
  $G$-module. Let  $p$ be a prime  dividing $|G|$ and  coprime to $q$,
  and $m$ an integer with $d/2 <m\leq d$.  An element $g\in G$ is said
  to be \emph{$(p,m)$-abundant}  if, in its action on  $V$, $g$ has an
  eigenvalue $\zeta$  in some extension field  of $\F_{q^\delta}$ such
  that $\zeta$ has multiplicative order divisible by $p$ and either
\begin{itemize}
\item[(i)]  $\zeta$ has $m$ Galois conjugates over $\F_{q^\delta}$, or
\item[(ii)]  $G\ne \GL_n(q)$, $m$ is even, $\zeta$ and $\zeta^{-1}$ are not Galois conjugate, 
and $\zeta$ and $\zeta^{-1}$ have together $m$ Galois conjugates over $\F_{q^\delta}$.
\end{itemize}
The element $g$ is  called \emph{$(p,m)$-abundant irreducible} in case
(i),  and \emph{$(p,m)$-abundant quasi-irreducible}  in case  (ii). In
either   case,  a   \emph{$p$-abundant}  element   is  one   which  is
$(p,m)$-abundant for some $m$ with $d/2<m\leq d$.  }
\end{Def}

The  terms  ``irreducible'' and  ``quasi-irreducible''  are chosen  to
reflect certain properties of  the actions of $p$-abundant elements on
the  natural $G$-module.   The  $(p,m)$-abundant irreducible  elements
leave invariant  a unique irreducible  subspace of dimension  $m$.   
In particular, we note that the $p$-abundant irreducible elements contain the family of so-called 
{\em primitive prime divisor} elements which underly the first and third 
authors' classical recognition algorithm~\cite{NPClassRec}. 
The  $(p,m)$-abundant
quasi-irreducible  elements  have  a  similar property,  preserving  a
specific decomposition of  a unique invariant $m$-dimensional subspace
into two closely related irreducible subspaces of dimension $m/2$.  The
proofs of these  facts are omitted here for brevity,  but may be found
in   our   related  paper~\cite{alg}   concerning
algorithms for identifying $p$-abundant elements computationally; 
see also the papers by Huppert~\cite{Hup1,Hup2}.

\begin{Thm} \label{Thm:intro}
Let $q$ be a prime power, $n$  an integer with $n \geq 9$, and $G$ as
in one of  the lines of Table~$\ref{tab:c2}$.  Suppose  that $p$ is an
odd  prime dividing  $|G|$  and coprime  to  $q$. Let  $e$ denote  the
smallest positive integer  such that $p$ divides $q^e-1$,  and $t$ the
largest    integer   such   that    $p^t$   divides    $q^e-1$.    Let
$Q(p;\textnormal{\texttt{I}};G)$  denote the  set of  all $p$-abundant
irreducible elements in $G$, and $Q(p;\textnormal{\texttt{QI}};G)$ the
set  of   all  $p$-abundant  quasi-irreducible   elements.   Then  for
$\textnormal{\texttt{T}}                     \in                    \{
\textnormal{\texttt{I}},\textnormal{\texttt{QI}}   \}$  and  constants
$c$,  $\kappa$, $\alpha$,  $\beta$  depending  on $G$  and  $e$ as  in
Table~$\ref{tab:c2}$, we have
\begin{equation} \label{eq:mainThmIneq}
-\left( \alpha + \frac{\beta\ln(2)}{p^t} \right) \frac{1}{n} - \frac{3\ln(2)}{e q^{n/4}} 
< \frac{|Q(p;\textnormal{\texttt{T}};G)|}{|G|} -  \left(1-\frac{1}{p^{t-1} (p+\kappa)}\right) \cdot \frac{c\ln(2)}{e} 
\leq \frac{\alpha}{n}.
\end{equation}
\end{Thm}

\begin{table}[t]
\begin{center}
\begin{tabular}{lllllllll}
\toprule
line & $G$ & \texttt{T} & $e$ parity & $e$ range & $c$ & $\kappa$ & $\alpha$ & $\beta$ \\
\midrule
\midrule
1 & \multirow{2}{*}{$\GL_n(q)$} & \multirow{2}{*}{\texttt{I}} & \multirow{2}{*}{all} & $e > n/p$ & $1$ & $0$ & $1$ & $0$ \\
2 & & & & $e \le n/p$ & $1$ & $1$ & $1$ & $1$ \\
\midrule
3 & \multirow{6}{*}{$\GU_n(q)$} & \multirow{2}{*}{\texttt{I}} & \multirow{2}{*}{$2\;\textnormal{(mod $4$)}$} & $e > 2n/p$ & $1$ & $0$ & $2$ & $0$ \\
4 & & & & $e \le 2n/p$ & $1$ & $1$ & $2$ & $1$ \\
\cline{3-9}
5 & & \multirow{4}{*}{\texttt{QI}} & \multirow{2}{*}{even} & $e > n/p$ & $1$ & $0$ & $1$ & $0$ \\
6 & & & & $e \le n/p$ & $1$ & $1$ & $1$ & $1$ \\
\cline{4-9}
7 & & & \multirow{2}{*}{odd} & $e > n/p$ & $1/2$ & $0$ & $1$ & $0$ \\
8 & & & & $e \le n/p$ & $1/2$ & $1$ & $1$ & $1/2$ \\
\midrule
9 & \multirow{4}{*}{$\begin{array}{l} \Sp_{2n}(q) \\ \SO_{2n+1}(q) \\ \SO_{2n}^\pm(q) \end{array}$} & \multirow{2}{*}{\texttt{I}} & \multirow{2}{*}{even} & $e > 2n/p$ & $1/2$ & $0$ & $3/2$ & $0$ \\
10 & & & & $e \le 2n/p$ & $1/2$ & $1$ & $3/2$ & $1/2$ \\
\cline{3-9}
11 & & \multirow{2}{*}{\texttt{QI}} & \multirow{2}{*}{all} & $e > n/p$ & $1/2$ & $0$ & $1$ & $0$ \\
12 & & & & $e \le n/p$ & $1/2$ & $1$ & $1$ & $1/2$ \\
\bottomrule
\end{tabular}
\end{center}
\caption{Cases for Theorem~\ref{Thm:intro}. (Line numbers are listed for later reference.)}\label{tab:c2}
\end{table}

The    proof    of    Theorem~\ref{Thm:intro}    is    given    in
Section~\ref{sec:proof},    using   preliminary theoretical  results
summarised  in  Section~\ref{sec:strat}  and technical lemmas collected in Section~\ref{sec:lemmas}. 
We note that we also take the opportunity to mention a small improvement to the results 
of our aforementioned paper~\cite{NPP} in Remark~\ref{rem:preinv} 
(that paper is unrelated to the present one, but also relies on 
the theory outlined in Section~\ref{sec:strat}).
Here  we just make a few remarks about  Theorem~\ref{Thm:intro}.

\begin{Rem}\label{rem:1}
{\rm  (a) For  combinations of  $\textnormal{\texttt{T}}$ and  $e$ not
  appearing        in       Table~$\ref{tab:c2}$,        the       set
  $Q(p;\textnormal{\texttt{T}};G)$        is        empty.        From
  Definition~\ref{def:pab-galois}, the quasi-irreducible case does not
  arise  in   $\GL_n(q)$.  For  $\GU_n(q)$,   $p$-abundant  irreducible
  elements arise only when $e \equiv 2\,(\textnormal{mod}\,4)$, and in
  the symplectic  and orthogonal  groups, the irreducible  case arises
  only  when $e$  is  even. The  details  are given  in  Section~\ref{sec:proof}.
  
  (b) A perhaps surprising consequence of Theorem~\ref{Thm:intro} is that the proportion 
  of $p$-abundant elements is at least a constant multiple of the lower bound 
  obtained by Isaacs et~al.~\cite{IKS95} for the proportion of {\em all} $p$-singular elements in $G$. 
  Specifically, upon observing that $1/(p^{t-1}(p+\kappa)) \leq 1/p$ in (\ref{eq:mainThmIneq}), we obtain
  \[
  \frac{|Q(p;\textnormal{\texttt{T}};G)|}{|G|} >  \left(1-\frac{1}{p}\right) \frac{c'}{e} 
  \]
  for some constant $c'$, in all cases where $Q(p;\textnormal{\texttt{T}};G)$ is nonempty.
  
(c) We do not consider the prime $p=2$, for which the results would be
  a  little  different  (Lemma~\ref{def:p|} would  need  modification,
  amongst  other things).   We  do, however,  believe  that a  similar
  result holds  in this case.
  
(d)    The    assumption    $n\geq    9$ is made  for technical
reasons,     as    certain     inequalities    used     in    deriving
(\ref{eq:mainThmIneq}) are  invalid for very small values  of $n$ (see
Lemmas~\ref{La:n>}  and \ref{La:partialProofs}(ii)).  However,  in proving
Theorem~\ref{Thm:intro} we obtain general  closed-form expressions
for    proportions     of    $p$-abundant    elements,     given    in
equations \eqref{eq:GLprf1}, \eqref{eq:GUprf1}, \eqref{eq:Spprf1odd} and \eqref{eq:Spprf1even}.  These expressions  depend  on certain  auxiliary
quantities  which  we estimate (using Lemma~\ref{la:SingerAll}) in order to  obtain  the  bounds given  in  the
theorem. But in  principle, they can be used  to calculate proportions
of $p$-abundant elements exactly, at least in certain simple cases. In
addition to  the small  $n$ cases not  covered by (\ref{eq:mainThmIneq}),  
we  have in  mind  situations  where $e$  is
reasonably  large,  say at  least  a  constant  fraction of  $n$.  
An example, where $G=\GL_n(q)$  with $e  \geq n/2$, is discussed in Remark~\ref{rem:GLexact} as illustration.
} 
\end{Rem}

\section{Strategy} \label{sec:strat}

Throughout the paper we use the following hypothesis.

\begin{hypothesis} \label{HypA}
\textnormal{
Let the group $G$, its dimension $d$ and the value of $\delta$ be 
as in one of the lines of Table~\ref{tab:G}, and let $V=V(d,q^{\delta})$ be the natural $G$-module. 
Assume that we have obtained $G$ as the fixed point 
set $\hat{G}^F$ of a connected reductive algebraic group $\hat{G}$ defined over the algebraic 
closure $\bar\F_q$ of $\F_q$, with $F$ a Frobenius morphism of $\hat{G}$. Moreover, assume 
that $F$ and a maximal torus $T_0$ in $\hat{G}$ have been chosen in the same way as outlined 
in \cite[Section 3]{LNP}, so that $W = N_{\hat{G}}(T_0)/T_0$ is the corresponding Weyl group.
}
\end{hypothesis}

\subsection{Estimating proportions via quokka sets}\label{ss:quokka}

In order to derive upper and lower bounds for proportions of $p$-abundant elements in the groups $G$ listed in Table~\ref{tab:G}, we apply the theory of {\em quokka sets} of finite groups of Lie 
type \cite{LNP,NPquokka}. 
These are subsets whose proportion in $G$ can be derived by determining certain proportions in maximal tori in $G$ and certain proportions in the corresponding Weyl group. 

Recall \cite[p.~11]{Carter85} that each element $g\in G$ has a unique Jordan decomposition $g =  su$, where $s\in G$ is semisimple, $u\in G$ is unipotent and $su=us$, with $s$ called the \emph{semisimple part} of $g$ and $u$ the \emph{unipotent part}. 
Note that the order $o(s)$ of $s$ is coprime to the characteristic, and $o(u)$ is a power of the characteristic.

The concept of a quokka set is introduced for finite groups of Lie type in \cite[Definition~1.1]{NPquokka}. 
A nonempty subset $Q$ of one of the groups $G$ in Table~\ref{tab:G} is a \emph{quokka set} if the following two conditions hold:
\begin{itemize}
\item[(i)] if $g\in G$ has Jordan decomposition $g=su$ with semisimple part $s$ and unipotent part $u$, then $g\in Q$ if and only if $s \in Q$;
\item[(ii)] $Q$ is a union of $G$-conjugacy classes.
\end{itemize}

We assume that Hypothesis~$\ref{HypA}$ holds, and summarise the required results. 
A subgroup $H$ of the connected reductive algebraic group $\hat{G}$ is said to be \emph{$F$-stable} if $F(H)=H$, and for each subgroup $H$ of $\hat{G}$ we write $H^F=H\cap G^F$. 
We define an equivalence relation on $W$ as follows: elements $w,w'\in W$ are \emph{$F$-conjugate} if there exists $x\in W$ such that  $w' =x^{-1} w F(x)$. 
The equivalence classes of this relation on $W$ are called \emph{$F$-conjugacy classes} \cite[p.~84]{Carter85}. 
The $G$-conjugacy classes of $F$-stable maximal tori are in one-to-one correspondence with the $F$-conjugacy classes of the Weyl group $W$ of $\hat{G}$. 
The explicit correspondence is given in \cite[Proposition 3.3.3]{Carter85}.

Let $\mathcal{C}$ be the set of $F$-conjugacy classes in $W$ and, for each $C\in\mathcal{C}$, let $T_C$ be a representative element of the family of $F$-stable maximal tori corresponding to $C$. 
The following theorem is a direct consequence of \cite[Theorem 1.3]{NPquokka}.

\begin{Thm}\label{the:quokka}
If Hypothesis~$\ref{HypA}$ holds and $Q \subseteq G$ is a quokka set, then
\[
\frac{|Q|}{|G|} =  \sum_{C \in \calC}\frac{|C|}{|W|} \frac{|T_C^F\cap  Q|}{|T_C^F|}.
\]
\end{Thm}

We refer to an $F$-stable maximal torus containing an element of $Q$ as a {\em quokka torus}, and call the corresponding $F$-conjugacy classes of $W$ {\em quokka classes}. 

In order to apply Theorem~\ref{the:quokka}, we check that the $p$-abundant elements in $G$ form quokka sets. 
We introduce the following notation, similar to that used in Theorem~\ref{Thm:intro} (for suitable $p, m$): $Q(p,m;\textnormal{\texttt{I}};G)$ denotes the set of all $(p,m)$-abundant irreducible elements in $G$, and $Q(p,m;\textnormal{\texttt{QI}};G)$ the set of all $(p,m)$-abundant quasi-irreducible elements. 
For brevity we combine this notation as $Q(p,m;\textnormal{\texttt{T}};G)$, where $\textnormal{\texttt{T}}$ is one of the symbols $\textnormal{\texttt{I}}$, $\textnormal{\texttt{QI}}$ (as in Theorem~\ref{Thm:intro}).
We have the following lemma.

\begin{La}\label{lem:quokkaset}
Suppose that Hypothesis~$\ref{HypA}$ holds. 
Let $p$ be a prime that is coprime to $q$, $\textnormal{\texttt{T}}\in \{ \textnormal{\texttt{I}},\textnormal{\texttt{QI}} \}$, and $m$ an integer with $d/2< m \leq d$ such that $Q(p,m;\textnormal{\texttt{T}};G)$ is nonempty. 
Then $Q(p,m;\textnormal{\texttt{T}};G)$ is a quokka set.
\end{La}

\begin{Prf}
It is clear from Definition~\ref{def:pab-galois} that each $Q(p,m;\textnormal{\texttt{T}};G)$ is a union of $G$-conjugacy classes. 
The condition that $g$ lies in $Q(p,m;\textnormal{\texttt{T}};G)$ if and only its semisimple part $s$ does follows from the well-known fact that $g$ and $s$ share the same characteristic polynomial.
\end{Prf}

Since the requirement that $m>d/2$ implies that an element can be $(p,m)$-abundant for at most one value of $m$, the sets $Q(p;\textnormal{\texttt{T}};G)$ in Theorem~\ref{Thm:intro} are then the disjoint unions of the respective $Q(p,m;\textnormal{\texttt{T}};G)$ over all $m$ with $d/2 < m \le d$.
Hence we have the following immediate corollary. 

\begin{Cor}\label{cor:quokkaset}
The nonempty $Q(p;\textnormal{\texttt{T}};G)$ are quokka sets and satisfy
\[
|Q(p;\textnormal{\texttt{T}};G)| = \sum_{d/2<m\leq d} |Q(p,m;\textnormal{\texttt{T}};G)|.
\]
\end{Cor}

\subsection{Maximal tori of the groups in Table {\bf \ref{tab:G}}} \label{subsec:tori}

Suppose that Hypothesis~\ref{HypA} holds. In order to apply Theorem~\ref{the:quokka} to estimate the proportion of elements in $G$ that lie in $Q(p;\textnormal{\texttt{T}};G)$ for some odd prime $p$ dividing $|G|$ and not dividing $q$, we have to describe the $F$-conjugacy classes of the Weyl group $W$ and their corresponding maximal tori. We summarise the description given in \cite{LNP} and \cite{NPquokka}, where more details can be found.

Consider first the cases where $G= \GL_n(q)$ or $G=\GU_n(q)$.  
Write $\delta=1$, $\epsilon =1$ in the first case and $\delta=2$, $\epsilon=-1$ in the second. Note that the Weyl group $W$ is isomorphic to $S_n$ via an isomorphism $W\rightarrow S_n$ that we denote (for reference in Section~\ref{ss:GL}) by $\sigma$. 
The $F$-conjugacy classes of $W$ are the conjugacy classes of $W$. So the $F$-conjugacy classes are parameterised by the partitions of $n$ describing the cycle types of permutations in $S_n$.
If $w\in W$ corresponds to the partition $\mu = (m_1,\ldots,m_k)$ of $n$ then each maximal torus $T^F$ of $G$ corresponding to the conjugacy class of $w$ is isomorphic to 
\[
\Z_{q^{m_1} - \epsilon^{m_1} }\times \cdots\times \Z_{q^{m_k} -\epsilon^{m_k}}.
\]
A cyclic direct factor $\Z_{q^m\pm 1}$ of $T^F$ corresponds to elements that have $m$ eigenvalues in $\bar{\F}_{q^\delta}$ that lie in $\F_{q^{\delta m}}$ and are permuted by the map $a \mapsto a^{\epsilon q}$. 
In particular, if $m>d/2$ and $T^F$ contains $(p,m)$-abundant elements, then for such elements these $m$ eigenvalues are precisely the Galois conjugates of $\zeta$ and $\zeta^{-1}$ as described in Definition~\ref{def:pab-galois}. 
(Note here that for $G=\GU_n(q)$, and similarly for the symplectic and orthogonal groups discussed below, if $\zeta$ is an eigenvalue of $g\in G$ then $\zeta^{-1}$ is also an eigenvalue, because $g$ is conjugate to its inverse transpose via the matrix of the form preserved by $G$.)

Now consider $G=\Sp_{2n}(q)$ or $G=\SO_{2n+1}(q)$. Here the Weyl group $W$ is isomorphic to $S_2 \wr S_{n}$, acting imprimitively on the set $\Omega=\{\pm1,\dots,\pm n\}$ of size $2n$, and consists of the so-called \emph{signed permutations}; that is to say, for $i,j\in\Omega$ and $g\in W$, $i^g=j$ if and only if $(-i)^g=-j$. 
We define a projection $\sigma : W\rightarrow S_n$ by mapping a signed permutation to the permutation it induces on $\{1, \ldots, n\}.$ 
For $g\in W$, a cycle of the image $\sigma(g)$ with length $\lambda$ is {\em positive} if it is the image under $\sigma$ of two $g$-cycles in $\Omega$ of length $\lambda$, and {\em negative} if it is the image under $\sigma$ of one $g$-cycle in $\Omega$ of length $2\lambda$.
A conjugacy class of $W$ is determined by its cycle type in $S_n$ and the numbers of positive cycles of each length. 
Suppose that $(\mu^+, \mu^-)$ is a partition of $n$ 
that determines a conjugacy class whose positive cycle lengths make up the parts of $\mu^+=(m^+_1,\dots,m^+_j)$ and negative cycle lengths make up the parts of $\mu^-=(m^-_1,\dots,m^-_k)$. 
Then each corresponding maximal torus $T^F$ of $G$ is isomorphic to 
\begin{equation}\label{eq:cl}
\left( \prod_{i=1}^j \Z_{q^{m^+_i}-1} \right)\times \left( \prod_{i=1}^k\Z_{q^{m^-_i}+1}\right).
\end{equation}
Here a cyclic factor $\Z_{q^\lambda\pm 1}$ corresponds to elements that have $m=2\lambda$ eigenvalues in $\bar{\F}_q$ that lie in $\F_{q^m}$ and are permuted by the map $a \mapsto a^q$. 
If $m>d/2$ and $T^F$ contains $(p,m)$-abundant elements, then for such elements these $m$ eigenvalues are precisely the Galois conjugates of $\zeta$ and $\zeta^{-1}$ as in Definition~\ref{def:pab-galois}.

Finally, consider $G=\SO^\pm_{2n}(q)$. 
We can view this group as a subgroup of $\SO_{2n+1}(q)$. 
The Weyl group $W$ has index $2$ in the Weyl group of  $\SO_{2n+1}(q)$, which we denote by $W_B$ below; namely, $W$ is the intersection of $W_B$ with the alternating group on $\Omega$. 
An element $w \in W_B$ lies in $W$ if and only if it has an even number of negative cycles. 
Moreover, we choose an element $w_n \in W_B$ such that $W_B=W\dot\cup Ww_n$, as described in \cite[Section~3.4]{LNP}.
The $F$-conjugacy classes of $W$ correspond to partitions $(\mu^+, \mu^-)$ of $n$ such that $\mu^-$ has an even number of parts in the case $G=\SO^+_{2n}(q)$ or an odd number of parts (and hence $|\mu^-|>0$) in the case $G=\SO^-_{2n}(q)$. 
The corresponding maximal tori are isomorphic to the groups in the product displayed in~(\ref{eq:cl}), and similar comments about $p$-abundant elements in these tori apply.

\begin{Rem} \label{rem:preinv}
\textnormal{
We take this opportunity to mention a small improvement to our paper \cite{NPP}, which is also based on the {\em quokka theory} outlined above. 
The lower bounds obtained in that paper for the proportions of so-called {\em pre-involutions} in finite classical groups can in fact be multiplied by $2$ in the cases $G=\SO_{2n}^\pm(q)$. 
Specifically, the `$1/4$' in the last line of \cite[Table 1]{NPP} may be replaced by `$1/2$', and the right-hand sides of the inequalities in \cite[Theorem 1.5(iii) and Corollary 1.6(ii)]{NPP} may be multiplied by $2$. 
In the arguments in \cite[Section 4.6]{NPP} for the cases $G=\SO_{2n}^\pm(q)$, the result \cite[Lemma 2.3]{NPP} should have been applied in conjunction with the fact (mentioned on \cite[p. 1025]{NPP}) that the Weyl group of type $D_n$ has index $2$ in the Weyl group of type $B_n$, as in the proof of \cite[Lemma 4.2]{NPquokka}; this would have yielded the additional factor of $2$ in our lower bounds.
}
\end{Rem}

\section{Proof of Theorem~\ref{Thm:intro}} \label{sec:proof}

We now prove Theorem~\ref{Thm:intro} using the strategy described in Section~\ref{sec:strat}. To aid the exposition, we refer in several places to technical results whose proofs are given in Section~\ref{sec:lemmas}. The following notation is used frequently: for a positive integer $k$ and a prime $r$, $(k)_r$ denotes the highest power of $r$ that divides $k$.

\subsection{$G=\GL_n(q)$} \label{ss:GL}

First suppose that $G= \GL_n(q)$. According to Definition~\ref{def:pab-galois}, only the irreducible case ($\textnormal{\texttt{T}}=\textnormal{\texttt{I}}$) arises. 
Here $C$ is a quokka class for $Q(p;\textnormal{\texttt{I}};G)$ if and only if $T_C^F$ contains a direct factor $A_m \cong \Z_{q^m-1}$ and $p$ divides $q^m-1$ for some $m$ with $n/2 < m \le n$. 
By Lemma~\ref{def:p|}(i), $p$ divides $q^m-1$ if and only if $e$ divides $m$, where $e$ is the order of $q$ modulo $p$. 
Thus, for an element $g$ of $C$, the image $\sigma(g) \in S_n$ has at least one cycle of length $m=be$, for some $b$, with $n/2 < m \le n$, namely $n/(2e)<b<m/e$. 
Note that any permutation in $S_n$ can have at most one such cycle and the proportion of elements in $S_n$ with a cycle of length $m$ is
\begin{equation}\label{SnProp}
\binom{n}{m}(m-1)! (n-m)! = \frac{1}{m}.
\end{equation}
For each $F$-conjugacy class $C$ with a cycle of length $m$, the maximal torus $T_C^F$ of $G$ corresponding to $C$ can be expressed as $T_C^F = A_m \times B$, where $B$ is a product of cyclic groups of the form $\Z_{q^{b_k}-1}$ with all $b_k\leq d/2$ and hence contains no $p$-abundant elements. 

Now, an element $g$ of this maximal torus $T_C^F$ is $(p,m)$-abundant (irreducible) if and only if, in its action on $V$, $g$ has an eigenvalue in some extension field of $\F_q$ with multiplicative order divisible by $p$ and $m$ Galois conjugates over $\F_q$. 
So, by the discussion in Section~\ref{subsec:tori}, $g$ is $(p,m)$-abundant irreducible if and only if its $A_m$-component has an eigenvalue in $\F_{q^m}$ with $m$ Galois conjugates over $\F_q$ and order divisible by $p$. 
Hence, denoting by $\theta(b)$ the proportion of elements in $\F^*_{q^{be}}$ that have $be$ Galois conjugates over $\F_q$ and order divisible by $p$, we have the following explicit expression for the proportion of $p$-abundant irreducible elements in $G=\GL_n(q)$:
\begin{equation} \label{eq:GLprf1}
\frac{|Q(p;\textnormal{\texttt{I}};G)|}{|G|} = \sum_{n/(2e) < b \leq n/e} \frac{\theta(b)}{be}.
\end{equation}

An upper bound for $\theta(b)$ is obtained by excluding the elements in $\F_{q^{be}}^*$ with order not divisible by $p$. 
These comprise the unique subgroup of $\F^*_{q^{be}}$ of index $|\F^*_{q^{be}}|_p = (q^{be}-1)_p$. 
By Lemma~\ref{def:p|}(iii), $(q^{be}-1)_p = p^{t+j}$, where $p^j=(b)_p$ and $p^t=(q^e-1)$, and hence $\theta(b) \leq 1-1/p^{t+j}$.
A lower bound for $\theta(b)$ is obtained by considering the proportion of elements in $\F_{q^{be}}^*$ with $be$ Galois conjugates, and then subtracting the proportion of elements in $\F_{q^{be}}^*$ with order not divisible by $p$.  
A lower bound for the former proportion is given by Lemma~\ref{la:SingerAll}(i) with $\ell=be$, yielding $\theta(b) > 1-1/p^{t+j}-3/q^{be/2}$.
Therefore, and since $be=m>d/2=n/2$, we have
\begin{equation} \label{la:torusPropGL}
1-\frac{1}{p^{t+j}}-\frac{3}{q^{n/4}} < \theta(b) \leq 1-\frac{1}{p^{t+j}}.
\end{equation}

We now estimate the sum in (\ref{eq:GLprf1}) to derive the bounds for $|Q(p;\textnormal{\texttt{I}};G)|/|G|$ claimed in Theorem~\ref{Thm:intro}. 
First suppose that $e>n/p$. 
Then all possible values of $b$ in (\ref{eq:GLprf1}) satisfy $b \leq n/e < p$ and hence have $j=0$ in the inequalities for $\theta(b)$ in \eqref{la:torusPropGL}. 
That is, we have $1-1/p^t-3/q^{n/4} < \theta(b) \leq 1-1/p^t$, independently of $b$. 
Then, using the notation and bounds of Lemma~\ref{la:propSumL}(i), we obtain
\begin{align}
\nonumber
\frac{|Q(p;\textnormal{\texttt{I}};G)|}{|G|} & \leq \left( 1-\frac{1}{p^t} \right) \frac{P(n/e,1)}{e} \\
\label{eq:GLprf2}
& \leq \left( 1-\frac{1}{p^t} \right) \left( \frac{\ln(2)}{e} + \frac{1}{n} \right) < \left( 1 - \frac{1}{p^t} \right) \frac{\ln(2)}{e} + \frac{1}{n}
\end{align}
and, similarly,
\begin{equation} \label{eq:GLprf3}
\frac{|Q(p;\textnormal{\texttt{I}};G)|}{|G|} > \left( 1 - \frac{1}{p^t} - \frac{3}{q^{n/4}} \right) \frac{P(n/e,1)}{e} > \left( 1 - \frac{1}{p^t} \right) \frac{\ln(2)}{e} - \frac{1}{n} - \frac{3\ln(2)}{eq^{n/4}}.
\end{equation}
Note that in \eqref{eq:GLprf3} we require that $1-1/p^t - 3/q^{n/4} > 0$, which holds under the assumption $n \geq 9$ made in Theorem~\ref{Thm:intro} according to Lemma~\ref{La:n>}.

Now consider the case where $e \leq n/p$. 
Denote by $i$ the positive integer satisfying $p^i \leq n/e < p^{i+1}$. 
Then each $b$ in (\ref{eq:GLprf1}) satisfies $(b)_p = p^j$ for some $j \in \{0,\ldots,i\}$, and hence the bounds for $\theta(b)$ given in \eqref{la:torusPropGL} depend on the variable $j=j(b)$ (unlike when $e>n/p$). 
We take this dependence on $j$ into account in order to obtain the precise leading term of $|Q(p;\texttt{I};G)|/|G|$. Write (\ref{eq:GLprf1}) as
\begin{equation} \label{genericSum}
\frac{|Q(p;\textnormal{\texttt{I}};G)|}{|G|} = \sum_{j=0}^i \left( \sum_{\substack{n/(2e) < b \leq n/e \\ (b)_p = p^j}} \frac{\theta(b)}{be} \right).
\end{equation}
Applying \eqref{la:torusPropGL} and using the notation of Lemma~\ref{la:propSumL}(i) yields
\begin{align}
\nonumber
\frac{|Q(p;\textnormal{\texttt{I}};G)|}{|G|} & \leq \frac{1}{e} \sum_{j=0}^i \left( \left( 1-\frac{1}{p^{t+j}} \right) \sum_{\substack{n/(2e) < b \leq n/e \\ (b)_p = p^j}} \frac{1}{b} \right) \\
\nonumber
& \leq \frac{1}{e} \sum_{j=0}^{i-1} \left( 1-\frac{1}{p^{t+j}} \right) \left( P(n/e,p^j) - P(n/e,p^{j+1}) \right) \\
\label{genericSumUB}
& \qquad {}+ \frac{1}{e} \left( 1-\frac{1}{p^{t+i}} \right) P(n/e,p^i)
\end{align}
and, similarly,
\begin{align}
\nonumber
\frac{|Q(p;\textnormal{\texttt{I}};G)|}{|G|} & > \frac{1}{e} \sum_{j=0}^{i-1} \left( 1-\frac{1}{p^{t+j}} - \frac{3}{q^{n/4}} \right) \left( P(n/e,p^j) - P(n/e,p^{j+1}) \right) \\
\label{genericSumLB}
& \qquad {}+ \frac{1}{e} \left( 1-\frac{1}{p^{t+i}} - \frac{3}{q^{n/4}} \right) P(n/e,p^i).
\end{align}
The bounds asserted in Theorem~\ref{Thm:intro} now follow upon application of Lemma \ref{la:propSumL}(i), and of Lemma~\ref{La:partialProofs} with $f_j = P(n/e,p^j)$ and hence $\ell=n/e$, $k_1 = \ln(2)$, $k_2=1$ (and $p,q,t,i$ as already defined). 
(Note that the assumption $n \geq 9$ made in the theorem is used when applying Lemma~\ref{La:partialProofs}(ii).)

\begin{Rem} \label{rem:GLexact}
\textnormal{
As mentioned in Remark~\ref{rem:1}(d), the closed-form expression \eqref{eq:GLprf1} for the proportion of $p$-abundant irreducible elements in $G=\GL_n(q)$ can, in principle, be used to compute this proportion exactly, at least in some simple cases. 
We have in mind situations where $e$ is reasonably large. As illustration, consider the case where $n/2<e\leq n$. 
The sum in \eqref{eq:GLprf1} then ranges over $b$ with $1/2<b<2$, so that $b=1$ is the only possible value, and one can check that $\theta(1)=1-1/p^t$ (by taking $\ell=1$ in Lemma~\ref{la:SingerAll}(i) and noting that the proportion considered there is then equal to $1$). 
For simplicity we do not use these facts in obtaining the estimates given in Theorem~\ref{Thm:intro}, but here they show that for $n/2<e\leq n$, the proportion of $p$-abundant irreducible elements in $\GL_n(q)$ is exactly $(1-1/p^t)/e$. 
Note that this particular case also follows from \cite[Lemma 5.6]{NPClassRec}. 
Similar comments apply to equations \eqref{eq:GUprf1}, \eqref{eq:Spprf1odd} and \eqref{eq:Spprf1even} below.
}
\end{Rem}

\subsection{$G=\GU_n(q)$} \label{ss:GU}

Now take $G=\GU_n(q)$. 
Here $C$ is a quokka class for $Q(p;\textnormal{\texttt{T}};G)$ if and only if $T_C^F$ has a direct factor $A_m \cong \Z_{q^m-(-1)^m}$ and $p$ divides $q^m-(-1)^m$ for some $m$ with $n/2<m\leq n$. 
So, depending on the parity of $m$, the image $\sigma(g)\in S_n$ of an element $g$ of $C$ must have a (unique) cycle of length $m$ as described below. 
In each case we obtain an expression analogous to \eqref{eq:GLprf1}.

If $m$ is odd then we need $p$ to divide $q^m+1$. By Lemma~\ref{def:p|}(ii), this occurs if and only if $e$ divides $2m$ and $e$ does not divide $m$. 
So $e \equiv 2\; \textnormal{(mod $4$)}$ and $m=be/2$ for some odd $b$ with $n/e<b\leq 2n/e$. 
The $(p,m)$-abundant elements in $T_C^F$ are of irreducible type, since condition (ii) of Definition~\ref{def:pab-galois} does not arise for $m$ odd. 
An element of $T_C^F$ is $(p,m)$-abundant irreducible if and only if its $A_m$-component has an eigenvalue in $\F_{q^{2m}}^*$ with $m$ Galois conjugates over $\F_{q^2}$ and multiplicative order divisible by $p$.
Hence, recalling from \eqref{SnProp} the proportion of elements in $S_n$ with an $m$-cycle, and denoting by $\theta_1(b)$ the proportion of elements in $\Z_{q^{be/2}+1} < \F_{q^{be}}^*$ with $be/2$ Galois conjugates over $\F_{q^2}$ and order divisible by $p$, we obtain the first case of \eqref{eq:GUprf1} below.

If $m$ is even then we need $p$ to divide $q^m-1$, which occurs if and only if $e$ divides $m$ (see Lemma~\ref{def:p|}(i)). 
Thus $m=be$ for some $b$ with $n/(2e)<b\leq n/e$, where $b$ must be even if $e$ is odd. 
The $(p,m)$-abundant elements in $T_C^F$ are of quasi-irreducible type, and an element of $T_C^F$ is $(p,m)$-abundant quasi-irreducible if and only if its $A_m$-component has an eigenvalue $\zeta$ in $\F_{q^m}^*$ with multiplicative order divisible by $p$ such that $\zeta$ and $\zeta^{-1}$ are not Galois conjugate and have together $m$ Galois conjugates over $\F_{q^2}$. 
Therefore, denoting by $\theta_2(b)$ the proportion of elements in $\F_{q^{be}}^*$ that have multiplicative order divisible by $p$, are not Galois conjugate to their inverses and have $be/2$ Galois conjugates over $\F_{q^2}$, we obtain the second and third cases below:
\begin{equation} \label{eq:GUprf1}
\frac{|Q(p;\textnormal{\texttt{T}};G)|}{|G|} = \left\{ \begin{array}{ll} 
\displaystyle{\sum_{\substack{n/e < b \leq 2n/e \\ \textnormal{$b$ odd}}} \frac{2\theta_1(b)}{be}} & \textnormal{if $e \equiv 2\; \textnormal{(mod $4$)}$ and $\textnormal{\texttt{T}}=\textnormal{\texttt{I}}$} \\
\displaystyle{\sum_{n/(2e) < b \leq n/e} \frac{\theta_2(b)}{be}} & \textnormal{if $e$ is even and $\textnormal{\texttt{T}}=\textnormal{\texttt{QI}}$} \\
\displaystyle{\sum_{\substack{n/(2e) < b \leq n/e \\ \textnormal{$b$ even}}} \frac{\theta_2(b)}{be}} & \textnormal{if $e$ is odd and $\textnormal{\texttt{T}}=\textnormal{\texttt{QI}}$}.
\end{array} \right.
\end{equation}

First consider the case where $\textnormal{\texttt{T}}=\textnormal{\texttt{QI}}$ with $e$ even.
Bounds on $\theta_2$ are obtained in a similar fashion to the bounds on $\theta$ in \eqref{la:torusPropGL}.
An upper bound $\theta_2(b) \leq 1 - 1/p^{t+j}$ is obtained by excluding the elements in $\F_{q^{be}}^*$ with order not divisible by $p$, which comprise the unique subgroup of index $(q^{be}-1)_p = p^{t+j}$, where $p^j=(b)_p$ and $p^t=(q^e-1)$.
A lower bound for $\theta_2(b)$ is obtained by considering the proportion of elements $\zeta \in \F_{q^{be}}^*$ that are not Galois conjugate to $\zeta^{-1}$ and have $be/2$ Galois conjugates over $\F_{q^2}$, and then subtracting the proportion of elements in $\F_{q^{be}}^*$ with order not divisible by $p$.  
A lower bound for the former proportion is given by Lemma~\ref{la:SingerAll}(ii) with $q$ replaced by $q^2$ and $\ell=be/2$, yielding $\theta_2(b) > 1-1/p^{t+j}-3/q^{be/2} > 1-1/p^{t+j}-3/q^{n/4}$. 
In other words, \eqref{la:torusPropGL} holds if $\theta$ is replaced by $\theta_2$.
Moreover, the corresponding (second) sum in (\ref{eq:GUprf1}) is the same as the sum in (\ref{eq:GLprf1}), except with $\theta$ replaced by $\theta_2$. 
It follows that we can obtain the same bounds for $|Q(p;\texttt{QI};G)|/|G|$ as we obtained for $|Q(p;\texttt{I};\GL_n(q))|/|\GL_n(q)|$. 
That is, the calculations in (\ref{eq:GLprf2})--(\ref{eq:GLprf3}) and (\ref{genericSum})--(\ref{genericSumLB}), with $\texttt{I}$ replaced by $\texttt{QI}$ and $\theta$ replaced by $\theta_2$ in \eqref{genericSum}, yield \eqref{eq:mainThmIneq} with the constants $c$, $\kappa$, $\alpha$, $\beta$ given in lines 5 and 6 of Table~\ref{tab:c2}, which are identical to lines 1 and 2, respectively.

Now let $\textnormal{\texttt{T}}=\textnormal{\texttt{QI}}$ with $e$ odd. The only difference from the previous case is that now the corresponding sum in (\ref{eq:GUprf1}) is restricted to even values of $b$. For $e>n/p$ this means that we proceed as in (\ref{eq:GLprf2})--(\ref{eq:GLprf3}) but with $\texttt{I}$ replaced by $\texttt{QI}$ and the proportion $P(n/e,1)$ replaced by $P'(n/e,1)$, where bounds on $P'$ are given in Lemma~\ref{la:propSumL}(ii). The result is that $\ln(2)$ is replaced by $\ln(2)/2$, and thus in line 7 of Table~\ref{tab:c2} as compared with line 5, the value of $c$ is divided by $2$. For $e \leq n/p$ we use (\ref{genericSum})--(\ref{genericSumLB}) with $\texttt{I}$ replaced by $\texttt{QI}$, $\theta$ replaced by $\theta_2$ in \eqref{genericSum}, $P(n/e,p^j)$ replaced by $P'(n/e,p^j)$ for $j=0,\ldots,i$, and the sums over $b$ restricted to even values of $b$. So now when applying Lemma~\ref{La:partialProofs} we set $f_j=P'(n/e,p^j)$ and hence $k_1=\ln(2)/2$ instead of $k_1=\ln(2)$. The result is that the values of both $c$ and $\beta$ are divided by $2$ in line 8 of Table~\ref{tab:c2} as compared with line 6.

It remains to consider the case where $\textnormal{\texttt{T}}=\textnormal{\texttt{I}}$, which arises (only) when $e \equiv 2\; \textnormal{(mod $4$)}$. 
The basic steps are similar to those in the preceding cases, with a few differences in the details.
An upper bound for $\theta_1(b)$ is obtained by excluding the elements of $\Z_{q^{be/2}+1}$ with order not divisible by $p$, which comprise the unique subgroup of $\Z_{q^{be/2}+1}$ of index $(q^{be/2}+1)_p$. 
By Lemma~\ref{def:p|}(iv), $(q^{be/2}+1)_p = p^{t+j}$, where $p^j=(b)_p$ and $p^t = (q^e-1)$, and so $\theta_1(b) \leq 1 - 1/p^{t+j}$. 
A lower bound for $\theta_1(b)$ is obtained by considering the proportion of elements $\Z_{q^{be/2}+1} < \F_{q^{be}}^*$ that have $be/2$ Galois conjugates over $\F_{q^2}$, and then subtracting the proportion of elements with order not divisible by $p$. 
A lower bound on former proportion is given by Lemma~\ref{la:SingerAll}(iv) with $\ell=be/2$, yielding $\theta_1(b) > 1 - 1/p^{t+j} - 3/q^{be/4} > 1 - 1/p^{t+j} - 3/q^{n/4}$, where the second inequality holds since $be=2m>n$ in the present case.
In summary, we have
\begin{equation} \label{la:torusPropU}
1-\frac{1}{p^{t+j}}-\frac{3}{q^{n/4}} < \theta_1(b) \leq 1-\frac{1}{p^{t+j}}.
\end{equation}

We now estimate the corresponding (first) sum in \eqref{eq:GUprf1}.
First suppose that $e>2n/p$. 
Then $b\leq 2n/e < p$ for all $b$ in the corresponding sum in (\ref{eq:GUprf1}), and so $j=0$ for all $b$ in the bounds for $\theta_1$ in \eqref{la:torusPropU}. 
Hence, in the notation of Lemma~\ref{la:propSumL}(iii), 
\[
\left(1-\frac{1}{p^t}-\frac{3}{q^{n/4}}\right) \frac{2P''(2n/e,1)}{e} < 
\frac{|Q(p;\textnormal{\texttt{I}};G)|}{|G|}
\leq \left(1-\frac{1}{p^t}\right) \frac{2P''(2n/e,1)}{e}.
\]
Applying Lemma~\ref{la:propSumL}(iii) and a calculation similar to (\ref{eq:GLprf2})--(\ref{eq:GLprf3}) yields \eqref{eq:mainThmIneq} with $c$, $\kappa$, $\alpha$, $\beta$ as in line 3 of Table~\ref{tab:c2}.

Now suppose that $e \leq 2n/p$. 
Let $i$ be the positive integer such that $p^i \leq n/e < p^{i+1}$. 
Each $b$ in the first sum in (\ref{eq:GUprf1}) satisfies $(b)_p = p^j$ for some $j \in \{0,\ldots,i\}$, and so the bounds in \eqref{la:torusPropU} depend on $j=j(b)$. 
We write
\[
\frac{|Q(p;\textnormal{\texttt{I}};G)|}{|G|} = \sum_{j=0}^i \left( \sum_{\substack{n/e < b \leq 2n/e \\ \textnormal{$(b)_p = p^j$, $b$ odd}}} \frac{2\theta_1(b)}{be} \right)
\]
and apply \eqref{la:torusPropU} to obtain the following inequalities analogous to (\ref{genericSumUB})--(\ref{genericSumLB}):
\begin{align}
\nonumber
\frac{|Q(p;\textnormal{\texttt{I}};G)|}{|G|} & \leq \frac{2}{e} \sum_{j=0}^{i-1} \left( 1-\frac{1}{p^{t+j}} \right) \left( P''(2n/e,p^j) - P''(2n/e,p^{j+1}) \right) \\
\nonumber
& \qquad {}+ \frac{2}{e} \left( 1-\frac{1}{p^{t+i}} \right) P''(2n/e,p^i), \\
\nonumber
\frac{|Q(p;\textnormal{\texttt{I}};G)|}{|G|} & > \frac{2}{e} \sum_{j=0}^{i-1} \left( 1-\frac{1}{p^{t+j}} - \frac{3}{q^{n/4}} \right) \left( P''(2n/e,p^j) - P''(2n/e,p^{j+1}) \right) \\
\nonumber
& \qquad {}+ \frac{2}{e} \left( 1-\frac{1}{p^{t+i}} - \frac{3}{q^{n/4}} \right) P''(2n/e,p^i).
\end{align}
Applying Lemma~\ref{la:propSumL}(iii), and Lemma~\ref{La:partialProofs} with $f_j = P''(2n/e,p^j)$ and hence $\ell=2n/e$, $k_1 = \ln(2)/2$, $k_2=2$ (and $p,q,t,i$ as already defined), we obtain \eqref{eq:mainThmIneq} with constants as in line 4 of Table~\ref{tab:c2}.

\subsection{$G=\Sp_{2n}(q)$, $\SO_{2n+1}(q)$ or $\SO_{2n}^\pm(q)$}

First suppose that $G=\Sp_{2n}(q)$ or $G=\SO_{2n+1}(q)$, and let $d=2n$ or $d=2n+1$, respectively. 
From the discussion in Section~\ref{subsec:tori}, here $C$ is a quokka class for $Q(p;\textnormal{\texttt{T}};G)$ if and only if $T_C^F$ has a direct factor $A_\lambda \cong Z_{q^\lambda\pm 1}$ and $p$ divides $q^\lambda\pm 1$ for some $\lambda$ such that $m=2\lambda$ satisfies $d/2<m\leq d$, which for the integer $\lambda$ is equivalent to $n/2<\lambda\leq n$.
The image $\sigma(g) \in S_n$ of an element $g$ of $C$ must have a cycle of length $\lambda$ as described below.

For a negative $\lambda$-cycle, $p$ must divide $q^{\lambda}+1$. 
According to Lemma \ref{def:p|}(ii), this occurs if and only if $e$ divides $2\lambda$ and $e$ does not divide $\lambda$. 
So $e$ must be even, and we need $\lambda = be/2$ for some odd $b$ with $n/2 < \lambda \leq n$, namely $n/e < b \leq 2n/e$. 
By \cite[Lemma~4.2(a)]{NPquokka}, the proportion of elements in $W$ with a negative cycle of length $\lambda$ is half the proportion of elements in $S_n$ with a cycle of length $\lambda$, and hence, by \eqref{SnProp}, is equal to $1/(2\lambda)$. 
The corresponding $(p,2\lambda)$-abundant elements in $T_C^F$ are of irreducible type, since a negative cycle of length $\lambda$ in $S_n$ corresponds to a single cycle of length $2\lambda$ in $W = S_2 \wr S_n$. 
An element of  $T_C^F$ is $(p,2\lambda)$-abundant irreducible if and only its $A_\lambda$-component has an eigenvalue in $\F_{q^{2\lambda}}$ with $2\lambda$ Galois conjugates over $\F_q$ and multiplicative order divisible by $p$. 
Denoting by $\theta^-(b)$ the proportion of elements in $\Z_{q^{be/2}+1} < \F_{q^{be}}^*$ with multiplicative order divisible by $p$ and $be$ Galois conjugates over $\F_q$, we obtain \eqref{eq:Spprf1even} below for the cases $G=\Sp_{2n}(q)$ and $G=\SO_{2n+1}$.

For a positive $\lambda$-cycle we need $p$ to divide $q^{\lambda}-1$, that is, $e$ must divide $\lambda$ (by Lemma~\ref{def:p|}(i)). 
So $\lambda = be$ for some $b$ with $n/2 < \lambda \leq n$, namely $n/(2e) < b \leq n/e$. 
By \cite[Lemma~4.2(a)]{NPquokka} and \eqref{SnProp}, the proportion of elements in $W$ with a positive cycle of length $\lambda$ is $1/(2\lambda)$. 
The corresponding $(p,2\lambda)$-abundant are of quasi-irreducible type, since a positive cycle of length $\lambda$ in $S_n$ corresponds to two cycles of length $\lambda$ in $W = S_2 \wr S_n$. 
An element of  $T_C^F$ is $(p,2\lambda)$-abundant quasi-irreducible if and only its $A_\lambda$-component has an eigenvalue $\zeta\in\F_{q^\lambda}$ with multiplicative order divisible by $p$ such that $\zeta$ and $\zeta^{-1}$ are not Galois conjugate and have together $2\lambda$ Galois conjugates over $\F_q$. 
Hence, denoting by $\theta^+(b)$ the proportion of elements $\zeta\in\F_{q^{be}}^*$ that have multiplicative order divisible by $p$, are not Galois conjugate to $\zeta^{-1}$ and have $be$ Galois conjugates over $\F_q$, we obtain \eqref{eq:Spprf1odd} below for the cases $G=\Sp_{2n}(q)$ and $G=\SO_{2n+1}$.

Now consider $G=\SO^\pm_{2n}(q)$, and recall the discussion at the end of Section~\ref{subsec:tori}.
A slight modification to the above argument is required. 
An $F$-conjugacy class $C$ in $W$ is a quokka class for $Q(p;\textnormal{\texttt{T}};G)$ if and only if, for an element $g$ of $C$, or of $Cw_n$ for $G=\SO^-_{2n}(q)$, the image $\sigma(g)\in S_n$ satisfies the same conditions as for $G=\Sp_{2n}(q)$ and in addition the total number of negative cycles is even if $G=\SO^+_{2n}(q)$ or odd if $G=\SO^-_{2n}(q)$. 
In particular, if $G=\SO^+_{2n}(q)$ then we cannot have a single negative cycle of length $n$, and for $G=\SO^-_{2n}(q)$ we cannot have a single positive cycle of length $n$. 
The result is that when $\lambda=n$, the Weyl group proportion $1/(2\lambda)=1/(2n)$ used above changes to $1/n$ if $\texttt{T}=\texttt{I}$, $G=\SO_{2n}^-(q)$ or $\texttt{T}=\texttt{QI}$, $G=\SO_{2n}^+(q)$ and to $0$ in the other two cases. 
This is reflected in equations \eqref{eq:Spprf1odd} and \eqref{eq:Spprf1even} below.

Let us summarise. 
Let $G=\Sp_{2n}(q),$ $\SO_{2n+1}(q)$ or $\SO_{2n}^\pm(q)$ for the remainder of this section. 
For the $p$-abundant quasi-irreducible elements in $G$, which arise for all values of $e$, we have
\begin{equation} \label{eq:Spprf1odd}
\frac{|Q(p;\textnormal{\texttt{QI}};G)|}{|G|} = \left( \sum_{n/(2e) < b \leq n/e} \frac{\theta^+(b)}{2be} \right) + \frac{\lambda^+\theta^+(n/e)}{2n},
\end{equation}
where
\[
\lambda^+ = \left\{ \begin{array}{rl}
1 & \textnormal{if $G=\SO^+_{2n}(q)$ and $e$ divides $n$} \\
-1 & \textnormal{if $G=\SO^-_{2n}(q)$ and $e$ divides $n$} \\
0 & \textnormal{otherwise.}
\end{array} \right.
\]
For even values of $e$ (only), we also have $p$-abundant irreducible elements arising, and the proportion of these elements is given by
\begin{equation} \label{eq:Spprf1even}
\frac{|Q(p;\textnormal{\texttt{I}};G)|}{|G|} = \left( \sum_{\substack{n/e < b \leq 2n/e \\ \textnormal{$b$ odd}}} \frac{\theta^-(b)}{be} \right) + \frac{\lambda^-\theta^-(2n/e)}{2n},
\end{equation}
with
\[
\lambda^- = \left\{ \begin{array}{rl}
-1 & \textnormal{if $G=\SO^+_{2n}(q)$ and $2n/e$ is an odd integer} \\
1 & \textnormal{if $G=\SO^-_{2n}(q)$ and $2n/e$ is an odd integer} \\
0 & \textnormal{otherwise.}
\end{array} \right.
\]

It remains to estimate the above sums in order to complete the proof of Theorem~\ref{Thm:intro}.
First consider (\ref{eq:Spprf1odd}). Begin by observing that
\[
\left| \frac{|Q(p;\textnormal{\texttt{QI}};G)|}{|G|} \quad - \sum_{n/(2e) < b \leq n/e} \frac{\theta^+(b)}{2be} \right| \leq \frac{1}{2n}.
\]
The sum in the above inequality is the same as the sum in \eqref{eq:GLprf1}, except for the factor of $1/2$ and the fact that $\theta$ has been replaced by $\theta^+$. 
Moreover, the inequalities in \eqref{la:torusPropGL} also hold if $\theta$ is replaced by $\theta^+$, as can be seen by applying part (ii) Lemma~\ref{la:SingerAll} with $\ell=be$ instead of part (i).
It follows that $|Q(p;\textnormal{\texttt{QI}};G)|/|G|$ satisfies the same bounds as those obtained in Section~\ref{ss:GL} for $|Q(p;\textnormal{\texttt{I}};\GL_n(q))|/|\GL_n(q)|$, except that we must divide by $2$ and then add $1/(2n)$ to the upper bound and subtract $1/(2n)$ from the lower bound. Thus, comparing lines 11 and 12 of Table~\ref{tab:c2} with lines 1 and 2, the values of $c$ and $\beta$ are halved, while halving $\alpha=1$ and then adding $1/2$ yields $\alpha=1$ again.

Now consider (\ref{eq:Spprf1even}), noting that
\[
\left| \frac{|Q(p;\textnormal{\texttt{I}};G)|}{|G|} \quad - \sum_{\substack{n/e < b \leq 2n/e \\ \textnormal{$b$ odd}}} \frac{\theta^-(b)}{be} \right| \leq \frac{1}{2n}.
\]
The sum above is the same as the first sum in (\ref{eq:GUprf1}), except divided by $2$ and with $\theta_1$ replaced by $\theta^-$. Moreover, the bounds for $\theta_1$ in \eqref{la:torusPropU} also hold with $\theta_1$ replaced by $\theta^-$, which is verified by applying part (iii) of Lemma~\ref{la:SingerAll} with $\ell=be/2$ instead of part (iv). It follows that $|Q(p;\textnormal{\texttt{I}};G)|/|G|$ satisfies the same bounds as those obtained in Section~\ref{ss:GU} for $|Q(p;\textnormal{\texttt{I}};\GU_n(q))|/|\GU_n(q)|$ in the case $e \equiv 2\; \textnormal{(mod $4$)}$, except that we must divide by $2$ and then add $1/(2n)$ to the upper bound and subtract $1/(2n)$ from the lower bound. In other words, comparing lines 9 and 10 of Table~\ref{tab:c2} with lines 3 and 4, the values of $c$ and $\beta$ are halved, while halving $\alpha=2$ and then adding $1/2$ yields $\alpha=3/2$.

\section{Technical results} \label{sec:lemmas}

Here we collect the various technical results used in Section~\ref{sec:proof}. 
Section~\ref{ss:torusProp} gives the results used to estimate the proportions $\theta$, $\theta_1$, $\theta_2$ and $\theta^\pm$ in equations \eqref{eq:GLprf1}, \eqref{eq:GUprf1}, \eqref{eq:Spprf1odd} and \eqref{eq:Spprf1even}. 
Section~\ref{ss:tech} collects the estimates of the various sums used to complete the proof of Theorem~\ref{Thm:intro}.

\subsection{Torus proportions} \label{ss:torusProp}

\begin{La} \label{def:p|}
Let $q$ be a prime power, $p$ a prime not dividing $q$, $e$ the smallest positive integer such that $p$ divides $q^e-1$, and write $p^t=(q^e-1)_p$. Then
\begin{itemize}
\item[\textnormal{(i)}] $p$ divides $q^m-1$ if and only if $e$ divides $m$,
\item[\textnormal{(ii)}] $p$ divides $q^m+1$ if and only if $e$ divides $2m$ and $e$ does not divide $m$,
\item[\textnormal{(iii)}] if $p\neq 2$ and $b$ is a positive integer with $(b)_p = p^j$ then $(q^{be}-1)_p = p^{t+j}$,
\item[\textnormal{(iv)}] if $p\neq 2$, $e$ is even and $b$ is an odd positive integer with $(b)_p = p^j$ then $(q^{be/2}+1)_p = p^{t+j}$.
\end{itemize}
\end{La}

\begin{Prf}
For proofs of (i) and (ii), see \cite[Lemma 4.5]{NPquokka}.

(iii) If $p$ does not divide $b$, namely $j=0$, then $(q^{be}-1)_p = p^t$ since
\[
\frac{q^{be}-1}{q^e-1} = 1 + q^e + \ldots + q^{e(b-1)} \equiv b\; \textnormal{(mod $p$)} 
\not \equiv 0\; \textnormal{(mod $p$)}.
\]
The proof is completed by induction on $j$. Claim 1 below says that without loss of generality we may assume that $b = p^j$. (Here we also note that $p$ does not divide $e$, since $e$ divides $p-1$ and is therefore coprime to $p$.) It then suffices to apply Claim 2 below with $k = q^{ep^j}$ and $r=t+j$.

{\em Claim $1$.} Let $v$, $v'$ be positive integers such that $v'$ is a multiple of $v$ and $(v')_p = (v)_p$. Then $(q^{v'}-1)_p = (q^v-1)_p$.

{\em Proof of Claim $1$.} Write $p^r = (q^v-1)_p$. If $v' = v\ell$ with $\ell$ not divisible by $p$ then
$(q^{v'}-1)/(q^v-1) = 1 + q^v + \ldots + q^{(\ell-1)v} \equiv \ell\; \textnormal{(mod $p^r$)}$, 
and therefore $(q^{v'}-1)/(q^v-1) \not \equiv 0\; \textnormal{(mod $p$)}.$

{\em Claim $2$.} If $k$ is a positive integer with $(k-1)_p = p^r$, where $r \geq 1$, then $(k^p-1)_p = p^{r+1}$.

{\em Proof of Claim $2$.} Write $k' = 1 + k + \ldots + k^{p-1}$. Then $k^p-1 = (k-1) k'$ and hence 
$(k^p-1)_p = p^r(k')_p$. So it suffices to check that $(k')_p = p$. Since $k = 1 + p^ry$ for some $y$ that is not divisible by $p$, we have
\begin{align*}
k' &=  1 + (1+yp^r) + (1+yp^r)^2 + \ldots + (1+yp^r)^{p-1} \\
&= p + yp^r(1+2+\ldots+(p-1)) + zp^{2r}
\end{align*}
for some $z$. So
\[
k' = p + y\frac{p-1}{2}p^{r+1} + zp^{2r},
\]
and since $r \geq 1$ it follows that $k'$ is divisible by $p$ but not by $p^2$, namely $(k')_p = p$.

(iv) Since $b$ is odd, $e$ does not divide $be/2$. So (i) implies that $p$ does not divide $q^{be/2}-1$, namely that $(q^{be/2}-1)_p = 1$, and then (iii) yields $(q^{be/2}+1)_p = (q^{be}-1)_p/(q^{be/2}-1)_p = p^{t+j}.$
\end{Prf}

\begin{La}\label{lem:GaloisCriterion}
Let $q$ be a prime power and $\ell$ a positive integer with $\ell \ge 2$. 
If $\zeta \in \F_{q^\ell}$ is Galois conjugate to $\zeta^{-1}$ over $\F_q$ then either $\zeta$
lies in a proper subfield of $\F_{q^\ell}$ or $\ell$ is even and
$\zeta$ lies in cyclic the subgroup of $\F_{q^{\ell}}^*$ of order $q^{\ell/2}+1$. 
\end{La}

\begin{Prf}
Let $i$ denote the least
nonnegative integer such that $\zeta^{q^i} = \zeta^{-1}$.
If $i=0$ then $\zeta^2 = 1$ and hence $\zeta=\pm 1$ lies in every proper
subfield of $\F_{q^\ell}$.  Suppose now that $i>0$. 
Then $\zeta^{q^{2i}} = \zeta$, or  $\zeta^{q^{2i}-1} = 1$. Since
  $\zeta^{q^{\ell}-1} = 1$, this implies that $\zeta^{q^{\gcd(\ell,2i)}-1}=
  \zeta^{\gcd(q^{\ell}-1,q^{2i}-1)}  =1$.  Write $k= \gcd(\ell,2i).$
If $k \not= \ell$ then $k<\ell$ and $k$ divides $\ell$, so $\zeta$ lies
in $\F_{q^k}<\F_{q^\ell}$ and the result
holds. This leaves
  $k = \ell$, in which case $i=r\ell/2$ for some positive integer $r$.
Observe that $\zeta^{q^{j\ell}}=\zeta$ for any nonnegative integer $j$.
If $r$ is even then
$\zeta^{-1} = \zeta^{q^i} = \zeta^{q^{\ell (r/2)}} =
\zeta$ and hence $i=0$, which we already considered above.
Suppose therefore that $r$ is odd, in which case $\ell$ is even.
Then $i = (r-1)/2 \cdot \ell + \ell/2$ and
$\zeta^{-1} = \zeta^{q^i} =  \zeta^{q^{(r-1)/2\cdot \ell+ \ell/2}} =
\zeta^{q^{\ell/2}}$. 
By the minimality of $i$, it follows that $r=1$ and so $i=\ell/2$. 
Thus $\zeta$ satisfies $\zeta^{q^{\ell/2}+1} =1$ and therefore lies in the cyclic
subgroup of $\F_{q^\ell}^*$ of order $q^{\ell/2}+1$.
\end{Prf}

\begin{La} \label{la:SingerAll}
Let $\ell$ be a positive integer. Then the following hold:
\begin{itemize}
\item[\textnormal{(i)}] The proportion
of elements in $\F_{q^\ell}^*$ with $\ell$ Galois conjugates
over $\F_q$ is greater than $1-3/q^{\ell/2}$.
\item[\textnormal{(ii)}]
The proportion of elements $\zeta \in \F_{q^\ell}^*$ that are not Galois conjugate 
to $\zeta^{-1}$ and have $\ell$ Galois conjugates
over $\F_q$ is greater than 
$1-3/q^{\ell/2}$. 
\item[\textnormal{(iii)}] The proportion of elements in
$\Z_{q^\ell+1} < \F_{q^{2\ell}}^*$ with $2\ell$ Galois
conjugates over $\F_{q}$ is greater than $1-3/
q^{\ell/2}$.
\item[\textnormal{(iv)}] For $\ell$ odd,  the proportion of elements in
$\Z_{q^\ell+1} < \F_{q^{2\ell}}^*$ with $\ell$ Galois
conjugates over $\F_{q^2}$ is greater than $1-3/
q^{\ell/2}$.
\end{itemize}
\end{La}

\begin{Prf}
(i) Write $A = \F_{q^\ell}^*$. The elements in $A$ with $\ell$
   Galois conjugates over $\F_q$ are precisely those that lie
   in no proper subfield of $A$.
Hence, denoting by $\rho(A)$ the set of all elements of $A$ that lie
in some field $K$ with $\F_q \leq K < \F_{q^\ell}^*$, we must show
that $|A\backslash \rho(A)|/|A| > 1-3/q^{\ell/2}$.
If $\ell=1$ then this inequality holds vacuously because $\rho(A)$ is
empty, so we now assume that $\ell \geq 2$.
If $\zeta$ is an element of some field $K$ with $\F_q \leq K <
\F_{q^\ell}^*$ then there is a prime divisor $r$ of $\ell$ such that
$\zeta \in K \leq \F_{q^{\ell/r}}^*$.
Hence
\[
|\rho(A)| \leq \sum_{\substack{\textnormal{$r$ prime} \\ r|\ell}} 
(q^{\ell/r}-1)
\leq \left( \sum_{j=1}^{\lfloor \ell/2 \rfloor} q^j \right) - 1 < 
\frac{q^{\ell/2+1}-1}{q-1} - 1 \leq 2q^{\ell/2}-2.
\]
So $|\rho(A)|/|A| < 2/q^{\ell/2}$, and thus $|A\backslash \rho(A)|/|A| > 1-2/q^{\ell/2}$.
In particular, $|A\backslash \rho(A)|/|A| > 1-3/q^{\ell/2}$ as claimed. 

(ii)
Write $A=\F_{q^\ell}^*$ again. If $\ell$ is
odd then by Lemma~\ref{lem:GaloisCriterion}  the proportion of
elements of $A$ that we are seeking is precisely as in (i), and is thus 
greater than $1-3/q^{\ell/2}$. Now suppose that $\ell$ is even.
Let $\rho(A)$ denote the set of elements in $A$ that lie in either a proper
subfield of $A$ or in the cyclic subgroup of $A$ of order
$q^{\ell/2}+1. $
By Lemma~\ref{lem:GaloisCriterion} the proportion that we are seeking is
$|A\backslash \rho(A)|/|A|$.  
From the proof of (i) we know that at most $2q^{\ell/2}-2$ elements lie in a proper subfield of $A$. 
So $|\rho(A)| \le (2 q^{\ell/2} - 2) + q^{\ell/2} = 3 q^{\ell/2}-2$, and 
thus $|A\backslash \rho(A)|/|A| > 1-3/q^{\ell/2}$.

(iii) 
Write $A = \Z_{q^\ell + 1}$ and let $\rho(A)$ denote the set of
all elements of $A$ that lie in some field $K$ with $\F_{q} \leq K <
\F_{q^{2\ell}}^*$. The proportion 
that we are seeking is $|A\backslash \rho(A)|/|A|$.  For
$\zeta\in \rho(A)$ we have $\zeta \in \F_{q^{2\ell/r}}^\ast$ for some
prime $r$ dividing $2\ell$, and so $|\zeta|$ divides
$$\operatorname{gcd}(q^\ell+1, q^{2\ell/r}-1) =
\begin{cases} q^{\ell/r}+1 &
\textrm{if } r \textrm{\ is odd} \\ 
\operatorname{gcd}(2,q-1) & \textrm{if } r = 2.
\end{cases}$$
If $\ell=1$  or $\ell$ is a power of $2$ then $r=2$ (only) and hence 
$|\rho(A)|/|A| \leq 2/(q^{\ell}+1)$. Then
$|A\backslash \rho(A)|/|A| \geq 1-2/q^{\ell} > 1-2/q^{\ell/2}$ and the 
result holds.
If $\ell$ is not a power of $2$ then
\begin{align*}
|\rho(A)| &\le  1 + \frac{2}{q^\ell+1} + \sum_{\substack{\textnormal{$r$ prime} \\ r|\ell, \;r\geq 3}}
q^{\ell/r}  \le \frac{2}{q^\ell+1} + \sum_{j=0}^{\lfloor \ell/3 \rfloor} q^j \\
&\le \frac{2}{q^\ell+1} + \frac{q^{\ell/3+1}-1}{q-1} \le  \frac{2}{q^\ell+1} + 2q^{\ell/3} < 3q^{\ell/3}.
\end{align*}
So $|A\backslash \rho(A)|/|A| > 1 - 3q^{\ell/3}/ q^{\ell} = 1 - 3/q^{2\ell/3} > 1 - 3/q^{\ell/2}$.

(iv) Write $A = \Z_{q^\ell + 1}$ and let $\rho(A)$ denote the set of
all elements of $A$ that lie in some field $K$ with $\F_{q^2} \leq K <
\F_{q^{2\ell}}^*$.
The result holds if $\ell=1$ since then $\rho(A)$ is empty, so we may 
assume that $\ell \geq 3$ (with $\ell$ odd). In this case an element in 
$\F_{q^{2\ell}}$ with $\ell$ Galois conjugates
over $\F_{q^2}$ has $2\ell$ Galois conjugates over $\F_{q}$, and so the 
result holds by (iii).
\end{Prf}

\subsection{Sums} \label{ss:tech}

\begin{La} \label{la:propSumL}
Let $\ell$ be a real number and $r$ an integer with $1 \leq r \leq \ell$. Write
\[
P(\ell,r) = \sum_{\substack{\ell/2 < b \leq \ell \\ r|b}} \frac{1}{b}, \quad 
P'(\ell,r) = \sum_{\substack{\ell/2 < b \leq \ell \\ \textnormal{$r|b$, $b$ even}}} \frac{1}{b}, \quad 
P''(\ell,r) = \sum_{\substack{\ell/2 < b \leq \ell \\ \textnormal{$r|b$, $b$ odd}}} \frac{1}{b},
\]
assuming further that $r$ is odd in the definition of $P''(\ell,r)$. Then
\begin{itemize}
\item[\textnormal{(i)}] $|P(\ell,r)-\ln(2)/r| \leq 1/\ell$,
\item[\textnormal{(ii)}] $|P'(\ell,r)-\ln(2)/(2r)| \leq 1/\ell$,
\item[\textnormal{(iii)}] $|P''(\ell,r)-\ln(2)/(2r)| \leq 2/\ell$.
\end{itemize}
\end{La}

\begin{Prf}
We make use of the following easily verified inequalities, in which $k_1$, $k_2$ are positive integers with $2 \leq k_1 \leq k_2$, and $x$ is a real number with $x>-1$: 
\begin{equation} \label{propSumLprf}
\ln\left(\frac{k_2+1}{k_1}\right) \leq \sum_{j=k_1}^{k_2}\frac{1}{j} \leq \ln\left(\frac{k_2}{k_1-1}\right), \quad \frac{x}{x+1} \leq \ln(1+x) \leq x.
\end{equation}
First consider the case where $r=1$. The above inequalities yield
\[
 \left\{ \begin{array}{ll}
-1/(\ell+1) & \textnormal{if $\ell$ is even} \\
0 & \textnormal{if $\ell$ is odd}
\end{array} \right. 
\leq P(\ell,1) - \ln(2)
\leq \left\{ \begin{array}{ll}
0 & \textnormal{if $\ell$ is even} \\
1/(\ell+1) & \textnormal{if $\ell$ is odd,}
\end{array} \right.
\]
where in the case of the upper bound for odd $\ell \geq 3$, 
we extract the $b=(\ell+1)/2$ term from the sum 
before applying the relevant inequality. 
It follows that $|P(\ell,1) - \ln(2)| \leq 1/(\ell+1) < 1/\ell$. 
If $r=2$ then writing $b=jr$ in the definition of $P(\ell,r)$ gives
\[
P(\ell,r) = \frac{1}{r} \sum_{\ell/(2r) < j \leq \ell/r} \frac{1}{j} = \frac{P(\lfloor \ell/r \rfloor,1)}{r},
\]
and it follows that $|P(\ell,r) - \ln(2)/r| \leq 1/(r(\lfloor \ell/r \rfloor+1)) \leq 1/\ell$. 
This completes the proof of (i). 
Assertions (ii) and (iii) follow: for $r=1$ we note that $P'(\ell,1) = P(\lfloor \ell/2 \rfloor,1)/2$ and $P''(\ell,1) = P(\ell,1) - P'(\ell,1)$, and then for $r=2$ we write $P'(\ell,r) = P'(\lfloor \ell/r \rfloor,1)/r$ and $P''(\ell,r) = P''(\lfloor \ell/r \rfloor,1)/r$.
\end{Prf}

\begin{La} \label{La:n>} 
Let $p$, $q$, $t$ be real numbers with $p \geq 3$, $q \geq 2$, $t \geq 1$. Then $1 - 1/p^t - 3/q^{n/4} > 0$ for all integers $n$ with $n\geq 9$.
\end{La}

\begin{Prf}
Since $p^t\geq 3$ and $q\geq 2$, the required inequality holds provided that $2/3 - 3/q^{n/4} > 0$, namely $n>4\log(9/2)/\log(2) \approx 8.68$.
\end{Prf}

\begin{La} \label{La:partialProofs} 
Let $p$ be a real number with $p \geq 3$ and $i$, $\ell$ positive integers with $p^i \leq \ell < p^{i+1}$. For $j\in\{0,\ldots,i\}$, suppose that $f_j$ are real numbers such that
\[
\frac{k_1}{p^j} - \frac{k_2}{\ell} \leq f_j \leq \frac{k_1}{p^j} + \frac{k_2}{\ell}
\]
for some positive real numbers $k_1,k_2$. Then
\begin{itemize}
\item[\textnormal{(i)}] for any real number $t$ with $t \geq 1$,
\[
\left( 1 - \frac{1}{p^{t+i}} \right) f_i + \sum_{j=0}^{i-1} \left( 1 - \frac{1}{p^{t+j}} \right) (f_j - f_{j+1}) < \left( 1 - \frac{1}{p^{t-1}(p+1)} \right) k_1 + \frac{k_2}{\ell};
\]
\item[\textnormal{(ii)}] if $n$ is an integer with $n\geq 9$ then for any real $q$, $t$ with $q \geq 2$ and $t \geq 1$,
\begin{align*}
&\left( 1 - \frac{1}{p^{t+i}} - \frac{3}{q^{n/4}} \right) f_i + \sum_{j=0}^{i-1} \left( 1 - \frac{1}{p^{t+j}} - \frac{3}{q^{n/4}} \right) (f_j - f_{j+1}) \\
>& \left( 1 - \frac{1}{p^{t-1}(p+1)} \right) k_1 - \left( k_2 + \frac{k_1}{p^t} \right) \frac{1}{\ell} - \frac{3k_1}{q^{n/4}}.
\end{align*}
\end{itemize}
\end{La}

\begin{Prf}
(i) We have
\begin{align*}
& \left( 1 - \frac{1}{p^{t+i}} \right) f_i + \sum_{j=0}^{i-1} \left( 1 - \frac{1}{p^{t+j}} \right) (f_j - f_{j+1}) = \left( 1 - \frac{1}{p^t} \right) f_0 + \frac{p-1}{p^t} \sum_{j=1}^{i} \frac{f_j}{p^j} \\
\leq& \left( 1-\frac{1}{p^t} \right) \left( k_1 + \frac{k_2}{\ell} \right) + \frac{p-1}{p^t} \sum_{j=1}^{i} \frac{1}{p^j} \left( \frac{k_1}{p^j} + \frac{k_2}{\ell} \right) \\
=& \left( 1-\frac{1}{p^t} + \frac{p-1}{p^t} \sum_{j=1}^{i} \frac{1}{p^{2j}} \right) k_1 + \left( 1-\frac{1}{p^t} + \frac{p-1}{p^t} \sum_{j=1}^{i} \frac{1}{p^j} \right) \frac{k_2}{\ell} \\
=& \left( 1 - \frac{1}{p^{t-1}(p+1)} - \frac{1}{p^t(p+1)p^{2i}} \right) k_1 + \left( 1-\frac{1}{p^tp^i} \right) \frac{k_2}{\ell} \\
<& \left( 1 - \frac{1}{p^{t-1}(p+1)} \right) k_1 + \frac{k_2}{\ell}.
\end{align*}
(ii) Set $E(n) = 3/q^{n/4}$ for brevity. By Lemma~\ref{La:n>}, the assumption $n \geq 9$ implies that $1 - 1/p^t - E(n) > 0$. This validates the first inequality below. The second inequality below follows from the assumption $p^{i+1} > \ell$, and the third from $p^i \leq \ell$, which implies that $\ell \geq p$ (since $i \geq 1$):
{\allowdisplaybreaks 
\begin{align*}
& \left( 1 - \frac{1}{p^{t+i}} - E(n) \right) f_i + \sum_{j=0}^{i-1} \left( 1 - \frac{1}{p^{t+j}} - E(n) \right) (f_j - f_{j+1}) \\
=& \left( 1 - \frac{1}{p^t} - E(n) \right) f_0 + \frac{p-1}{p^t} \sum_{j=1}^{i} \frac{f_j}{p^j} \\
>& \left( 1-\frac{1}{p^t} - E(n) \right) \left( k_1 - \frac{k_2}{\ell} \right) + \frac{p-1}{p^t} \sum_{j=1}^{i} \frac{1}{p^j} \left( \frac{k_1}{p^j} - \frac{k_2}{\ell} \right) \\
=& \left( 1-\frac{1}{p^t} + \frac{p-1}{p^t} \sum_{j=1}^{i} \frac{1}{p^{2j}} - E(n) \right) k_1 - \left( 1-\frac{1}{p^t} + \frac{p-1}{p^t} \sum_{j=1}^{i} \frac{1}{p^j} - E(n) \right) \frac{k_2}{\ell} \\
=& \left( 1-\frac{1}{p^{t-1}(p+1)} - \frac{1}{p^t(p+1)p^{2i}} - E(n) \right) k_1 - \left( 1-\frac{1}{p^tp^i} - E(n) \right) \frac{k_2}{\ell} \\
>& \left( 1 - \frac{1}{p^{t-1}(p+1)} - \frac{1}{p^{t-2}(p+1)\ell^2} - E(n) \right) k_1 - \frac{k_2}{\ell} \\
=& \left( 1 - \frac{1}{p^{t-1}(p+1)} \right) k_1 - \frac{k_2}{\ell} - \frac{k_1}{p^{t-2}(p+1)\ell^2} - k_1E(n) \\
>& \left( 1 - \frac{1}{p^{t-1}(p+1)} \right) k_1 - \left( k_2 + \frac{k_1}{p^t} \right) \frac{1}{\ell} - k_1E(n).
\end{align*}
}
\end{Prf}

\section*{Acknowledgements}

The authors are grateful to Simon Guest for pointing out an error in one of the lemmas in Section~\ref{sec:lemmas} in an early version of the paper, and to an anonymous referee for numerous invaluable comments that significantly improved the paper.
The first and third authors acknowledge the support of the Australian Research Council Discovery Projects DP0879134 and DP110101153. The second author acknowledges support within the Australian Research Council Federation Fellowship FF0776186 of the third author.

\end{document}